# AUTOMORPHISMS OF THE SEMIGROUP OF ENDOMORPHISMS OF FREE ASSOCIATIVE ALGEBRAS

A. KANEL-BELOV, A. BERZINS AND R. LIPYANSKI

ABSTRACT. Let $A = A(x_1, ..., x_n)$ be a free associative algebra in the variety of associative algebras $\mathcal{A}$ freely generated over $K$ by a set $X = \{x_1, ..., x_n\}$, $End\, A$ be the semigroup of endomorphisms of $A$, and $Aut\, End\, A$ be the group of automorphisms of the semigroup $End\, A$. We investigate the structure of the groups $Aut\, End\, A$ and $Aut\, \mathcal{A}^\circ$, where $\mathcal{A}^\circ$ is the category of finitely generated free algebras from $\mathcal{A}$. We prove that the group $Aut\, End\, A$ is generated by semi-inner and mirror automorphisms of $End\, F$ and the group $Aut\, \mathcal{A}^\circ$ is generated by semi-inner and mirror automorphisms of the category $\mathcal{A}^\circ$. This result solves an open Problem formulated in [14].

## 1. INTRODUCTION

Let $\Theta$ be a variety of linear algebras over a commutative-associative ring $K$ and $F = F(X)$ be a free algebra from $\Theta$ generated by a finite set $X$. Here $X$ is supposed to be a subset of some infinite universum $X^0$. Let $G$ be an algebra from $\Theta$ and $K_\Theta(G)$ be the category of algebraic sets over $G$. Here and bellow we refer to [15, 16] for Universal Algebraic Geometry (UAG) definitions used in our work.

The category $K_\Theta(G)$ can be considered from the point of view of the possibility to solve systems of equations in the algebra $G$. Algebras $G_1$ and $G_2$ from $\Theta$ are categorically equivalent if the categories $K_\Theta(G_1)$ and $K_\Theta(G_2)$ are correctly isomorphic. Algebras $G_1$ and $G_2$ are geometrically equivalent if

$$T''_{G_1} = T''_{G_2}$$

holds for all finite sets $X$ and for all binary relations $T$ on $F$ and $'$ is Galois correspondence between sets in $Hom(F, G)$ and the binary relations on $F$.

It has been shown in [16] that categorical and geometrical equivalences of algebras are related and their relation is determined by the structure of the group $Aut\, \Theta^0$, where $\Theta^0$ is the category of free finitely generated algebras of $\Theta$. There is a natural connection between a structure of the groups $Aut\, End\, F$, $F \in \Theta$, and $Aut\, \Theta^0$.

Let $\mathcal{A}$ be the variety of associative algebras with (or without) 1, $A = A(x_1, ..., x_n)$ be a free associative algebra in $\mathcal{A}$ freely generated over $K$ by a set $X = \{x_1, ..., x_n\}$. One of our aim here is to describe the group $Aut\, End\, A$ and, as a consequence, to obtain a description of the group $Aut\, \mathcal{A}^\circ$ for the variety of associative algebras over a field $K$.

We prove that the group $Aut\,End\,A$ is generated by semi-inner and mirror automorphisms of $End\,A$ and the group $Aut\,\mathcal{A}^\circ$ is generated by semi-inner and mirror automorphisms of the category $\mathcal{A}^\circ$.

Earlier, the description of $Aut\,\mathcal{A}^\circ$ for the variety $\mathcal{A}$ of associative algebras over algebraically closed fields has been given in [11] and, over infinite fields, in [3]. Also in the same works, the description of $Aut\,End\,F(x_1, x_2)$ has been obtained.

Note that a description of the groups $Aut\,End\,F$, $F \in \Theta$, and $Aut\,\Theta^\circ$ for some other varieties $\Theta$ has been given in [2, 3, 4, 6, 7, 8, 9, 11, 12, 13, 14, 18].

## 2. Automorphisms of the semigroup $End\,F$ and of the category $\Theta^0$

We recall the basic definitions we use in the case of the variety $\mathcal{A}$ of associative algebras over a field $K$.

Let $F = F(x_1, ..., x_n)$ be a finitely generated free algebra of a variety $\Theta$ of linear algebras over $K$ generated by a set $X = \{x_1, ..., x_n\}$.

**Definition 2.1.** [2] An automorphism $\Phi$ of the semigroup $End\,F$ of endomorphisms of $F$ is called quasi-inner if there exists a bijection $s : F \to F$ such that $\Phi(\nu) = s\nu s^{-1}$, for any $\nu \in End\,F$; $s$ is called adjoint to $\Phi$.

**Definition 2.2.** [15] A quasi-inner automorphism $\Phi$ of $End\,F$ is called semi-inner if its adjoint bijection $s : F \to F$ satisfies the following conditions:
  1. $s(a + b) = s(a) + s(b)$,
  2. $s(a \cdot b) = s(a) \cdot s(b)$,
  3. $s(\alpha a) = \varphi(\alpha)s(a)$,

for all $\alpha \in K$ and $a, b \in F$ and an automorphism $\varphi : K \to K$. If $\varphi$ is the identity automorphism of $K$, we say that $\Phi$ is an inner.

Let $A = A(x_1, ..., x_n)$ be a finitely generated free associative algebra over a field $K$ of the variety $\mathcal{A}$. Further, without loss of generality, we assume that associative algebras of $\mathcal{A}$ contain 1.

**Definition 2.3.** [11] A quasi-inner automorphism $\Phi$ of $End\,A$ is called mirror if its adjoint bijection $s : A \to A$ is anti-automorphism of $A$.

Recall the notions of category isomorphism and equivalence [10]. An isomorphism $\varphi : \mathcal{C} \to \mathcal{D}$ of categories is a functor $\varphi$ from $\mathcal{C}$ to $\mathcal{D}$ which is a bijection both on objects and morphisms. In other words, there exists a functor $\psi : \mathcal{D} \to \mathcal{C}$ such that $\psi\varphi = 1_{\mathcal{C}}$ and $\varphi\psi = 1_{\mathcal{D}}$.

Let $\varphi_1$ and $\varphi_2$ be two functors from $\mathcal{C}_1$ to $\mathcal{C}_2$. A functor isomorphism $s : \varphi_1 \longrightarrow \varphi_2$ is a collection of isomorphisms $s_D : \varphi_1(D) \longrightarrow \varphi_2(D)$ defined for all $D \in Ob\,\mathcal{C}_1$ such that for every $\nu : D \longrightarrow B$, $\nu \in Mor\,\mathcal{C}_1$, $B \in Ob\,\mathcal{C}_1$, holds

$$s_B \cdot \varphi_1(\nu) = \varphi_2(\nu) \cdot s_D,$$

i.e., the following diagram is commutative

$$\begin{array}{ccc} \varphi_1(D) & \xrightarrow{s_D} & \varphi_2(D) \\ \varphi_1(\nu) \downarrow & & \downarrow \varphi_2(\nu) \\ \varphi_1(B) & \xrightarrow{s_B} & \varphi_2(B) \end{array}$$

The isomorphism of functors $\varphi_1$ and $\varphi_2$ is denoted by $\varphi_1 \cong \varphi_2$.

An equivalence between categories $\mathcal{C}$ and $\mathcal{D}$ is a pair of functors $\varphi : \mathcal{C} \to \mathcal{D}$ and $\psi : \mathcal{D} \to \mathcal{C}$ together with natural isomorphisms $\psi\varphi \cong 1_{\mathcal{C}}$ and $\varphi\psi \cong 1_{\mathcal{D}}$. If $\mathcal{C} = \mathcal{D}$, then we get the notions of automiorphism and autoequivalence of the category $\mathcal{C}$.

For every small category $\mathcal{C}$ denote the group of all its automorphisms by $Aut\,\mathcal{C}$. We will distinguish the following classes of automorphisms of $\mathcal{C}$.

**Definition 2.4.** [7, 13] An automorphism $\varphi : \mathcal{C} \to \mathcal{C}$ is equinumerous if $\varphi(D) \cong D$ for any object $D \in Ob\,\mathcal{C}$ ; $\varphi$ is stable if $\varphi(D) = D$ for any object $D \in Ob\,\mathcal{C}$ ; and $\varphi$ is inner if $\varphi$ and $1_{\mathcal{C}}$ are naturally isomorphic, i.e., $\varphi \cong 1_{\mathcal{C}}$.

In other words, an automorphism $\varphi$ is inner if for all $D \in Ob\,\mathcal{C}$ there exists an isomorphism $s_D : D \to \varphi(D)$ such that

$$\varphi(\nu) = s_B \nu s_D^{-1} : \varphi(D) \to \varphi(B)$$

for any morphism $\nu : D \to B$, $B \in Ob\,\mathcal{C}$.

Let $\Theta$ be a variety of linear algebras over $K$. Denote by $\Theta^0$ the full subcategory of finitely generated free algebras $F(X), |X| < \infty$, of the variety $\Theta$.

**Definition 2.5.** [13] Let $A_1$ and $A_2$ be algebras from $\Theta$, $\delta$ be an automorphism of $K$ and $\varphi : A_1 \to A_2$ be a ring homomorphism of these algebras. A pair $(\delta, \varphi)$ is called semimomorphism from $A_1$ to $A_2$ if

$$\varphi(\alpha \cdot u) = \alpha^\delta \cdot \varphi(u), \ \ \forall \alpha \in K, \forall u \in A_1.$$

Define the notion of a semi-inner automorphism of the category $\Theta^0$.

**Definition 2.6.** [13] An automorphism $\varphi \in Aut\,\Theta^0$ is called semi-inner if there exists a family of semi-isomorphisms $\{s_{F(X)} = (\delta, \tilde{\varphi}) : F(X) \to \tilde{\varphi}(F(X)),\ F(X) \in Ob\,\Theta^0\}$, where $\delta \in Aut\,K$ and $\tilde{\varphi}$ is a ring isomorphism from $F(X)$ to $\tilde{\varphi}(F(X))$ such that for any homomorphism $\nu : F(X) \longrightarrow F(Y)$ the following diagram

$$\begin{array}{ccc} F(X) & \stackrel{s_{F(X)}}{\longrightarrow} & \tilde{\varphi}(F(X)) \\ \nu \downarrow & & \downarrow \varphi(\nu) \\ F(Y) & \underset{s_{F(Y)}}{\longrightarrow} & \tilde{\varphi}(F(Y)) \end{array}$$

is commutative.

Now we define the notion of a mirror automorphism of the category $\mathcal{A}^\circ$.

**Definition 2.7.** [16] An automorphism $\varphi \in Aut\,\mathcal{A}^\circ$ is called mirror if it does not change objects of $\mathcal{A}^\circ$ and for every $\nu : A(X) \to A(Y)$, where $A(X), A(Y) \in Ob\,\mathcal{A}^\circ$, it holds

$$\varphi(\nu) : A(X) \to A(Y) \ \text{ such that } \ \varphi(\nu)(x) = \delta(\nu(x)), \ \forall x \in X,$$

where $\delta : A(Y) \to A(Y)$ is the mirror automorphism of $A(Y)$.

Further, we will need the following

**Proposition 2.8.** [7, 13] For any equinumerous automorphism $\varphi \in Aut\,\mathcal{C}$ there exists a stable automorphism $\varphi_S$ and an inner automorphism $\varphi_I$ of the category $\mathcal{C}$ such that $\varphi = \varphi_S \varphi_I$.

## 3. Quasi-inner automorphisms of the semigroup $End\,F$ for associative and Lie varieties

We will need the standard endomorphisms of free algebra $F = F(x_1, ..., x_n)$ of the variety $\Theta$.

**Definition 3.1.** [9] Standard endomorphisms of $F$ in the base $X = \{x_1, ..., x_n\}$ are the endomorphisms $e_{ij}$ of $F$ which are determined on the free generators $x_k \in X$ by the rule: $e_{ij}(x_k) = \delta_{jk} x_i$, $x_i \in X$, $i, j, k \in [1n]$, $\delta_{jk}$ is the Kronecker delta.

Denote by $S_0$ a subsemigroup of $End\,F$ generated by $e_{ij}$, $i, j \in [1n]$. Further, we will use the following statements

**Proposition 3.2.** [9] Let $\Phi \in Aut\,End\,F(X)$. Elements of the semigroup $\Phi(S_0)$ are standard endomorphisms in some base $U = \{u_1, ..., u_n\}$ of $F$ if and only if $\Phi$ is a quasi-inner automorphism of $End\,F$.

The description of quasi-inner automorphisms of $End\,A(X)$, where $A(X)$ is a free associative algebra with or without 1 over a field K, is following

**Proposition 3.3.** [4, 9] Let $\Phi \in Aut\,End\,A(X)$ be a quasi-inner automorphism of $End\,A(X)$. Then $\Phi$ is either a semi-inner or a mirror automorphism, or a composition of them.

Let us investigate the images of standard endomorphisms under automorphisms of $End\,A$. To this end we introduce endomorphisms of rank 1.

**Definition 3.4.** We say that an endomorphism $\varphi$ of $A$ has rank 1, and write this as $rk\,(\varphi) = 1$, if its image $Im\,\varphi$ is a commutative subalgebra of $F$.

Note that according to Bergman's theorem [1], the centralizer of any non-scalar element of $A$ is a polynomial ring in one variable over $K$. Thus, an endomorphism $\varphi$ of $A$ is of rank 1 if and only if $\varphi(A) = K[z]$ for some element $z \in A$.

**Proposition 3.5.** An endomorphism $\varphi$ of $A$ is of rank 1 if and only if there exists a non-zero endomorphism $\psi \in End\,A$ such that for any $h \in End\,A$

$$\varphi \circ h \circ \psi = 0 \tag{3.1}$$

*Proof.* Let $\varphi \in End\,A$ be an endomorphism of rank 1. Let us take the endomorphism $\psi \in End\,A$ such that

$$\psi(x_1) = [x_1, x_2] \text{ and } \psi(x_i) = 0 \text{ for all } i \neq 1.$$

Since $\varphi(A)$ is a commutative subalgebra of $F$, the condition (3.1) is fulfilled for any $h \in End\,A$.

Conversely, let the condition (3.1) is fulfilled for the endomorphism $\varphi$. Assume, on the contrary, that $Im\,\varphi$ is not a commutative algebra. Without loss of generality, it can be supposed that $[\varphi(x_1), \varphi(x_2)] \neq 0$, $x_1, x_2 \in X$. Denote by $R = K[\varphi(x_1), \varphi(x_2)]$ a subalgebra of $A$ generated by $\varphi(x_1)$ and $\varphi(x_2)$. It is well known (see [5]) that $R$ is a free non-commutative subalgebra of $A$.

Since $\psi$ is a non-zero endomorphism of $A$, there exists $x_i \in X$ such that $\psi(x_i) \neq 0$. Set $P = P(x_1, ..., x_n) = \psi(x_i)$. We wish to show that $P$ is an identity of the algebra $R$. Assume, on the contrary, that there exist elements $z_1, ..., z_n \in R$ such that $P(z_1, ..., z_n) \neq 0$. Consider sets $\varphi^{-1}(z_i)$, $i \in [1n]$, and choose elements $y_i \in \varphi^{-1}(z_i), i \in [1n]$, from them. We may construct an endomorphism $h$ of $A$ such that $h(x_i) = y_i$, $i \in [1n]$. Then we have

$$0 = \varphi \circ h \circ \psi(x_i) = P(\varphi \circ h(x_1), ..., \varphi \circ h(x_n)) = P(\varphi(y_1), ..., \varphi(y_n)) = P(z_1, ..., z_n).$$

We arrived at a contradiction. Therefore, $P$ is an identity of $R$. Since $R$ is a free non-commutative subalgebra of $A$, it has no non-trivial identities. Thus, $P = 0$. We get a contradiction again. Therefore, $Im\,\varphi$ is a commutative algebra and Proposition is proved. □

It follows directly from this Proposition

**Corollary 3.6.** Let $\Phi \in Aut\,End\,A$ and $rk\,(\varphi) = 1$. Then $rk\,(\Phi(\varphi)) = 1$.

**Definition 3.7.** A set of endomorphisms $\mathcal{B}_e = \{e'_{ij} \mid e'_{ij} \in End\,A, i,j \in [1n]\}$ of $A$ is called a subbase of $End\,A$ if
  1. $e'_{ij} e'_{km} = \delta_{jk} e'_{im}$, $\forall i, j, k, m \in [1n]$,
  2. $rk\,(e'_{ij}) = 1$, $\forall i, j \in [1n]$, i.e., there exist elements $z_{ij} \in A$, $i, j \in [1n]$, such that $e'_{ij}(A(X)) = K[z_{ij}]$ for all $i, j \in [1n]$.

Further, for simplicity, we write $z_{ii} = z_i$, $i \in [1n]$.

**Definition 3.8.** We say that a subbase $\mathcal{B}_e$ is a base collection of endomorphisms of $A$ (or a base of $End\,A$, for short) if $Z = \langle z_i \mid z_i \in A, i \in [1n]\rangle$ is a base of $A$.

**Proposition 3.9.** A subbase of endomorphisms $\mathcal{B}_e$ is a base if and only if for any collection of endomorphisms $\alpha_i : A \to A$, $\forall i \in [1n]$, and any subbase $\mathcal{B}_f = \{f'_{ij} \mid i, j \in [1n]\}$ of $End\,A$ there exist endomorphisms $\varphi, \psi \in End\,A$ such that

(3.2) $$\alpha_i \circ f'_{ii} = \psi \circ e'_{ii} \circ \varphi, \text{ for all } i \in [1n].$$

*Proof.* Let a subbase of endomorphisms $\mathcal{B}_e$ be base. Since $rk\,(f'_{ij}) = 1$, $\forall i, j \in [1n]$, there exist elements $y_{ij} \in A$, $i, j \in [1n]$, such that $f'_{ij}(A(X)) = K[y_{ij}]$ for all $i, j \in [1n]$. We define an endomorphisms $\psi$ and $\varphi$ of $A$ in the following way:

$$\varphi(x_i) = z_i \text{ and } \psi(z_i) = \alpha_i(y_i), \text{ for all } i \in [1n]$$

where $y_i = y_{ii}$, $\forall i \in [1n]$. Since $Z = \langle z_i \mid z_i \in A, i \in [1n]\rangle$ is a base of $A$, the definition of the endomorphism $\psi$ is correct. Now, it is easy to check that the condition (3.2) with the given $\phi$ and $\psi$ is fulfilled.

Conversely, assume that the condition (3.2) is fulfilled for the subbase $\mathcal{B}_e$. Let us prove that $Z = \langle z_i \mid z_i \in A, i \in [1n]\rangle$ is a base of $A$. Choosing in (3.2) $\alpha_i = e_{ii}$ and $f'_{ij} = e_{ij}$ for all $i, j \in [1n]$, we obtain

$$e_{ii} = \psi \circ e'_{ii} \circ \varphi,$$

i.e., $\psi(e'_{ii}\varphi(x_i)) = x_i$ for all $i \in [1n]$. Denote $t_i = e'_{ii}\varphi(x_i)$. We have $\psi(t_i) = x_i$. Since $A$ is Hopfian, the elements $t_i$, $i \in [1n]$, form a base of $A$. Taking into account the equality $e'_{ii}(A(X)) = K[z_i]$, we obtain $t_i = \chi_i(z_i) \in K[z_i]$. Since $t_i z_i = z_i t_i$, $i \in [1n]$, by Bergman's theorem we have $z_i = g_i(t_i)$. Thus, $z_i = g_i(\chi_i(z_i))$. Similarly, $t_i = \chi_i(g_i(t_i))$. Therefore, there exists non-zero elements $a_i$ and $b_i$ in $K$ such that $z_i = a_i t_i + b_i$, $i \in [1n]$. Thus, $Z = \langle z_i \mid z_i \in F, \forall i \in [1n]\rangle$ is also a base of $A$ as claimed. □

Now we deduce

**Corollary 3.10.** Let $\Phi \in Aut\,End\,A$. Then $\mathcal{C} = \{\Phi(e_{ij}) \mid i, j \in [1n]\}$ forms a base collection of endomorphisms of $A$.

*Proof.* Since $\Phi(e_{ij})\Phi(e_{km}) = \delta_{jk}\Phi(e_{im})$ and by Corollary 3.6, $rk\,(\Phi(e_{ij})) = 1$, the set $\mathcal{C}$ is a subbase of $End\,A$. It is evident that the condition (3.2) is fulfilled for the subbase $\mathcal{C}$. By Proposition 3.9, $\mathcal{C}$ is a base of $End\,A$. $\square$

**Lemma 3.11.** *Let $\mathcal{B}_e = \{e'_{ij} \mid e'_{ij} \in End\,A,\, i,j \in [1n]\}$ be a base collection of endomorphisms of $End\,A$. Then there exists a base $S = \langle s_k \mid s_k \in A, k \in [1n] \rangle$ of $A$ such that the endomorphisms from $\mathcal{B}_e$ are standard endomorphisms in $S$.*

*Proof.* Since $(e'_{ii})^2 = e'_{ii}$, we have $e'_{ii}(z_i) = z_i$, $i \in [1n]$. The equality $e'_{ii}e'_{ij}z_j = e'_{ij}z_j$ implies the existence of a polynomial $f_j(z_i) \in K[z_i]$ such that $e'_{ij}z_j = f_j(z_i)$. Similarly, there exists a polynomial $g_i(z_j) \in K[z_j]$ such that $e'_{ji}z_i = g_i(z_j)$. We have
$$z_j = e'_{jj}z_j = e'_{ji}e'_{ij}z_j = e'_{ji}(f_j(z_i)) = f_j(g_i(z_j)) \text{ for all } i,j \in [1n].$$
and, in similar way, $z_i = g_i(f_j(z_i))$ for all $i,j \in [1n]$. Thus $f_j$ and $g_i$ are linear polynomials over $K$ in variables $z_i$ and $z_j$, respectively. Therefore,

(3.3) $$e'_{ij}z_j = a_j z_i + b_j,\ a_i, b_i \in K \text{ and } a_i \neq 0.$$

Note that $e'_{ij}z_k = e'_{ij}e'_{kk}z_k = 0$ if $k \neq j$. Now we have for $i \neq j$
$$0 = {e'_{ij}}^2 z_j = e'_{ij}(a_j z_i + b_j) = e'_{ij}((b_j)) = b_j,$$
i.e., $e'_{ij}z_j = a_j z_i$, $a_j \neq 0$. Let $V = Span(z_1, ..., z_n)$. Then $V$ is the vector space over $K$ with a basis $Z = \langle z_k \mid z_k \in A, k \in [1n] \rangle$ and $e'_{ij}$, $i,j \in [1n]$, are linear operators on $V$. Set
$$S = \langle s_i = e'_{i1}z_1 \mid z_1 \in Z,\, i \geq 1 \rangle.$$
Since $s_i = a_1 z_i$, $a_1 \neq 0$, $i \in [1n]$, we have that $S$ is a base of $A$. In this base we obtain $e'_{ij}s_k = \delta_{jk}s_i$, $i,j,k \in [1n]$. The proof is complete. $\square$

## 4. Structure of automorphisms of the semigroup $End\,F$ for associative and Lie varieties

Now we give the description of the groups $Aut\,End\,A$ and $Aut\,\mathcal{A}^{\circ}$.

**Theorem 4.1.** *The group $Aut\,End\,A$ is generated by semi-inner and mirror automorphisms of $End\,A$.*

*Proof.* By Corollary 3.10, the set of endomorphisms $\mathcal{C} = \{\Phi(e_{ij}) \mid \forall i \in [1n]\}$ is a base collection of endomorphisms of $A$. By Lemma 3.11, there exists a base $S = \langle s_k \mid s_k \in A, k \in [1n] \rangle$ such that the endomorphisms $\Phi(e_{ij})$ are standard endomorphisms in $S$. According to Proposition 3.2, we obtain that $\Phi$ is quasi-inner. By virtue of Proposition 3.3, the group $Aut\,End\,A$ is generated by semi-inner and mirror automorphisms of $End\,A$ as claimed. $\square$

Using Theorem 4.1 we prove

**Theorem 4.2.** *The group $Aut\,\mathcal{A}^{\circ}$ of automorphisms of the category $\mathcal{A}^{\circ}$ is generated by semi-inner and mirror automorphisms of the category $\mathcal{A}^{\circ}$.*

*Proof.* Let $\varphi \in Aut\,\mathcal{A}^{\circ}$. It is clear that $\varphi$ is an equinumerous automorphism. By Proposition 2.8, $\varphi$ can be represented as the composition of a stable automorphism $\varphi_S$ and an inner automorphism $\varphi_I$. Since a stable automorphism does not change free algebras from $\mathcal{A}^{\circ}$, we obtain that $\varphi_S \in Aut\,End\,A(x_1, ..., x_n)$. By Theorem 4.1, $\varphi_S$ is generated by semi-inner and mirror automorphisms of $End\,A$. Using this fact and Reduction Theorem [7, 13], we obtain that the group $Aut\,\mathcal{A}^\circ$ generated by semi-inner and mirror automorphisms of the category $\mathcal{A}^\circ$. This ends the proof. □

## 5. Acknowledgments

The authors are grateful to B. Plotkin for attracting their attention to this problem and interest to this work.

## References

[1] G. Bergman, Centralizers in free associative algebras, *Trans. Amer. Math. Soc.* **137** (1969) 327-344.
[2] A. Berzins, B. Plotkin, E. Plotkin, Algebric geometry in varieties of algebras with the given algebra of constants, *Journal of Math. Sciences* **102** (3) (2000) 4039-4070.
[3] A. Berzins, The group of automorphisms of the category of free associative algebra, *Preprint* (2004).
[4] A. Berzins, The group of automorphisms of semigroup of endomorphisms of free commutative and free associative algebra, *Preprint*, (2004)
[5] P. Cohn, *Free rings and their relations*, (Academic Press, London, 1985).
[6] E. Formanek, A question of B. Plotkin about the semigroup of endomorpjsms of a free group, *Proc. American Math. Soc.* **130** (2001) 935-937.
[7] Y. Katsov, R. Lipyanski, B. Plotkin, Automorphisms of categories of free modules and free Lie algebras , (2004) pp. 18, to appear.
[8] R. Lipyanski, B. Plotkin, Automorphisms of categories of free modules and free Lie algebras, *Preprint. Arxiv:math. RA//0502212 (2005)*.
[9] R. Lipyanski, Automorphisms of the semigroup of endomorphisms of free algebras of homogeneous varieties, *Preprint. Atxiv: math. RA//0511654v1* (2005).
[10] S. Mac Lane, *Categories for the Working Mathematician*, (New York-Berlin: Spinger-Verlag, 1971).
[11] G. Mashevitzky, Automorphisms of the semigroup of endomorphisms of free ring and free associative algebras, Preprint
[12] G. Mashevitzky, B. Schein, Automorphisms of the endomorphism semigroup of a free monoid or a free semigroup, *Proc. Amer. math. Soc.* **8** (2002) 1-10.
[13] G. Mashevitzky, B. Plotkin, E. Plotkin, Automorphisms of the category of free Lie algebras, *Journal of Algebra* **282** (2004) 490-512.
[14] G. Mashevitzky, B. Plotkin, E. Plotkin, Automorphisms of the category of free algebras of varieties, *Electron. Res. Announs. Amer. Math. Soc.* **8** (2002) 1-10.
[15] B. Plotkin, Seven lectures in universal algebraic geometry, *Preprint. Arxiv:math. RA/0502212* (2002).
[16] B. Plotkin, Algebra with the same (algebraic geometry), in *Proc. of the Steklov Institut of Mathematics* **242** (2003) 176-207.
[17] B. Plotkin, G. Zhitomirskii, On automorphisms of categories of free algebras of some varieties. *Preprint. Arxiv:math. RA/0501331* (2005).
[18] G. Zhitomirskii, Automorphisms of the semigroup of all endomorphisms of free algebras, *Preprint. Arxiv:math. GM/0510230 v1* (2005).

*Department of Mathematics*, *The Hebrew University, Jerusalem, 91904, Israel*
*E-mail address*: `kanel@huji.ac.il`

*University of Latvia, Riga LV 1586, Latvia*
*E-mail address*: `aberzins@latnet.lv`

*Department of Mathematics*, *Ben Gurion University, Beer Sheva, 84105, Israel*
*E-mail address*: `lipyansk@cs.bgu.ac.il`